\documentclass[12 pt]{report}
\usepackage[brazil]{babel}
\usepackage[latin1]{inputenc}
\usepackage{amsmath,amssymb,indentfirst}

\begin{document}
\centerline{{\bf AN ALTERNATIVE METHOD FOR THE UNDETERMINED COEFFICIENTS}}
\centerline{{\bf AND THE ANNIHILATOR METHODS}}

\

\centerline{\bf{Oswaldo Rio Branco de Oliveira}}

\vspace{0,3 cm}

\hspace{-0,6 cm}{\bf Abstract.} This paper exhibits a very simple formula for a particular solution of a linear ordinary differential equation with constant real coefficients, $P(\frac{d}{dt})x = f$, $f$ a function given by a linear combination of polynomials, trigonometrical and exponential real functions products. This is done by solving lower triangular linear systems with constant coefficients. Two examples are given.

\vspace{0,2 cm}

\hspace{- 0,6 cm}{\sl Mathematics Subject Classification: 34A30}

\hspace{- 0,6 cm}{\sl Key words and phrases:} Ordinary Differential Equations, Linear Equations and Systems.

\vspace{0,2 cm}

\centerline{\bf{Introduction}}

\vspace{0,1 cm}

The objective of this paper is to present a very simple and alternative method to the method of undetermined coefficients (see Birkhoff and Rota [2], Braun [3], Coddington and Levinson [4], Krantz and Simmons [5], Pontryagin [7], Robinson [8], Ross [9] and [10], and Zill [11]) and also to the annihilator method (see [10], [11] and Apostol [1]), both very well known, of solving a linear ordinary  differential equation with constant real coefficients, $P(\frac{d}{dt})x = f$, where $f$ is a real function given by a linear combination of polynomials, trigonometrical functions and exponential functions products. This is accomplished by reducing the given equation to the trivial case in which $f$ is a polynomial. Most text books develop the method of undetermined coefficients but not the annihilator method, considered to be a bit more sophisticated. In most books, it is not proven that the undetermined coefficients strategy is a valid one and such strategy is, in some texts, named ``the guessing method'' or ``the lucky guess method''. It is worth to point out that a paper by Ross, see [9], explains very nicely why the undetermined coefficients strategy works. In general, applying either the undetermined coefficients or the annihilator method  requires a good computations amount. There are three advantages in the alternative presented here. First, at many times it reduces the computations amount necessary in comparison with either the undetermined coefficients or the annihilator method. Second, the proof is quite easy. Third, it shows a very clear algorithm to find a particular solution for the given linear differential equation.  

\vspace{0,2 cm}

As it is well known, the general solution of the homogeneous equation associated to the given ode is a consequence of the Fundamental Theorem of Algebra (see Oliveira [6]), and can be found in the textbooks mentioned above. 

\

\centerline{{\bf Preliminaries}}

\vspace{0,3 cm}

Given an infinitely differentiable function $f:\mathbb R \to \mathbb C$, then $f^{(k)}$ indicates the $k$-th derivative of $f$. We also write, $f=f^{(0)}$, $f'=f^{(1)}$, $f''=f^{(2)}$, $f'''=f^{(3)}$, etc. Moreover, we indicate $C^{\infty}(\mathbb R;\mathbb C)=\{f:\mathbb R \to \mathbb C: f^{(k)}\ \textrm{exists}\,,\forall k \in \mathbb N\}$.

\vspace{0,2 cm}

\hspace{-0,6 cm}{\bf Lemma 1.} Let us consider the linear ordinary differential inhomogeneous equation with coefficients $a_j\in \mathbb C$, $0\leq j\leq n$, with at least one $a_j\neq 0$, $1\leq j\leq n$, 
\[a_nx^{(n)} + a_{n-1}x^{(n-1)} +\ ...\ +\ a_1x' + a_0x \,=\, R\ ,\,\ x = x(t)\,, \,t \in \mathbb R\,,\] 
where the function $R=R(t)=b_nt^n +...+b_1t +b_0$ is a non zero polynomial in the real variable $t$, with coefficients $b_j \in \mathbb C$, $0\leq j\leq n$. Then, the equation has
a polynomial solution $Q=Q(t): \mathbb R \to \mathbb C$. Moreover,
\begin{itemize}
\item[(a)] If $a_0 \neq 0$, then we have degree$(Q)=$ degree$(R)$.
\item[(b)] If $k =\max\{j: a_l= 0\,,l\leq j\}$ exists, then we can admit $Q = t^{k+1}Q_1$, where $Q_1$ is a polynomial satisfying degree$(Q_1) =$ degree$(R)$. 
\item[(c)] If all the coefficients above are real, then $Q$ and $Q_1$ are real polynomials. 
\end{itemize}
{\bf  Proof.} If $k$ exists, then it is clear that $k<n$.
\begin{itemize}
\item[(a)] Let us solve the pair of equations, 
\begin{displaymath}
\left\{\begin{array}{ll}
\ \ \ \ \ \ \ \ \ Q(t) = c_nt^n + ...+ c_1t + c_0,\\
(1.1)\ \ \ a_0Q \, +\,  a_1Q' \, +\, ...\ +\, a_jQ^{(j)}\ +\ ... \, +\,  a_nQ^{(n)} = R,
\end{array}
\right.
\end{displaymath}
by identifying the coefficient of $t^{n-l}$, $l\leq n$, at the parcels $a_jQ^{(j)}$, $j\leq l$, noticing that at the others parcels the coefficient of $t^{n-l}$ is zero. Fixed  the parcel $a_jQ^{(j)}$, $j \leq l$, a factor of the  coefficient comes up with the computation,
\newpage 
\[c_{n-l+j}\frac{d^j}{dt^j}\{\,t^{n-l+j}\} = c_{n-l+j}(n-l+j)(n-l+j -1)...(n-l+1)t^{n-l},\] 
and the coefficient looked for is then $a_jc_{n-l+j}\frac{(n-l+j)!}{(n-l)!}$.

Therefore, substituting in (1.1) the expression $Q(t)=c_nt^n + ...+ c_1t + c_0$, and its derivatives, we conclude that the coefficient of $t^{n-l}$ satisfy the identity (in which, the sum is in decreasing order on $j$, where $j=l,l-1,...,0$), 
\[(1.2)\ \ \ \ a_lc_n\frac{n!}{(n-l)!} \ +\ ...\ +\ a_jc_{n-l+j}\frac{(n-l+j)!}{(n-l)!}\ +\ ...\ +\  a_0c_{n-l} = b_{n-l}\ \,,\ l = 0,1,...,n \ .\]
From the expressions given by $(1.2)$ it follows the trivial matricial equation,
\begin{displaymath}
\left[\begin{array}{ccccccccccc}
a_0 & 0 & 0 &0 &0 &0 & 0 & . & 0 & 0 & 0\\
a_1n & a_0 & 0 &0 & 0 & 0 & 0 & .& 0 & 0 & 0 \\
a_2\frac{n!}{(n-2)!} & a_1\frac{(n-1)!}{(n-2)!} & a_0&0 & 0 & 0 & 0 & .& 0 & 0 & 0 \\
.& . & . &. &. &. & . & .& . & . & . \\
 .& . & . &. & . & . & . & .& . & . & . \\
 .& . & . &. & . & . & . & .& . & . & . \\
a_l\frac{n!}{(n-l)!} & a_{l-1}\frac{(n-1)!}{(n-l)!} &. & a_j\frac{(n-l+j)!}{(n-l)!} & . & .& a_0 & 0 & 0 & 0 & 0\\
. & . & . &. & . & . & . & .& . & . & . \\
. & . & . &. &. &. & . & .& a_0 & 0 & 0 \\
 . & . & . &. & . & . & . & .& . & a_0 & 0 \\
a_nn! & .& . & .& . & . & . & .& a_22! & a_1 & a_0 \\
\end{array}\right] 
\left[\begin{array}{l}
c_n\\
c_{n-1}\\
c_{n-2}\\
.\\
. \\
.\\
c_{n-l}\\
. \\
c_2 \\
c_1\\
c_0\\
\end{array}\right]
= \left[\begin{array}{l}
b_n\\
b_{n-1}\\
b_{n-2}\\
.\\
. \\
.\\
b_{n-l} \\
.\\
b_2\\
b_1\\
b_0
\end{array}\right].
\end{displaymath}
\item[(b)] The equation is $a_nx^{(n)} \,+\,...\,+\,a_{k+1}x^{(k+1)} = R$. By (a), the equation 
 $a_ny^{(n-k-1)} +...+ a_{k+1}y = R$ has a solution $y(t)=Q(t)$, with degree$(Q)$ = degree$(R)$. Integrating $(k+1)$-times the function $y = y(t)$ and choosing, at each time, zero for independent term, we reach the desired polynomial.
\item[(c)] It is trivial $\blacksquare$ 
\end{itemize}

\hspace{-0,6 cm}{\bf Lemma 2.} Let $P\big(\frac{d}{dt}\big)= a_n\frac{d^n}{dt^n}\ +\ a_{n-1}\frac{d^{n-1}}{dt^{n-1}}\, + ...+ \,a_1\frac{d}{dt} \ +\ a_0I$, where $a_j\in \mathbb R$, $0\leq j \leq n$, and $I$ is the identity operator over $C^{\infty}(\mathbb R;\mathbb C)$. Let $p(t)$ be the characteristic polynomial of $P\big(\frac{d}{dt}\big)$. Then, if $Q=Q(t) \in C^{\infty}(\mathbb R;\mathbb C)$ and $\gamma \in \mathbb C$, 
\[ P\big(\frac{d}{dt}\big)\left\{\,Q(t)e^{\gamma t}\,\right\}\ =\ \left[\frac{p^{(n)}(\gamma)}{n!}Q^{(n)} \, +\,  ...\, +\, \frac{p''(\gamma)}{2!}Q''\, +\,   \frac{p'(\gamma)}{1!}Q' \, +\,  \frac{p(\gamma)}{0!}Q \, \right]e^{\gamma t} \ . \]
{\bf Proof.} It is clear that
\[\big(\frac{d}{dt} - \gamma I\big)\{\,e^{\gamma t}Q(t)\,\} = e^{\gamma t }Q'(t)\,,\]
and iterating we arrive at
\[\big(\frac{d}{dt} - \gamma I\big)^k\{\,e^{\gamma t}Q(t)\,\} = e^{\gamma t}Q^{(k)}(t)\ .\]

Applying Taylor's formula to $p(t)$ about $\gamma$ we find,
\[ p(t) = \sum_{k=0}^n \frac{p^{(k)}(\gamma)}{k!}(t - \gamma)^k\,.\]
Substituting $\frac{d}{dt}$ for $t$ and computing the resulting operator at $e^{\gamma t}Q(t)$ gives
\[ P\big(\frac{d}{dt}\big)\{\,e^{\gamma t}Q(t)\,\} = \sum_{k=0}^n \frac{p^{(k)}(\gamma)}{k!}\big(\frac{d}{dt} -\gamma I\big)^k\{\,e^{\gamma t }Q(t)\,\} = \sum_{k=0}^n \frac{p^{(k)}(\gamma)}{k!}e^{\gamma t }Q^{(k)}(t)\  \blacksquare\]

\vspace{0,2 cm}

Let us keep the notation of Lemma 2. Let us consider $R$, a real non zero polynomial, and a pair of numbers $(\gamma,\delta)$ such that 
\[\bullet\ \ \ \ \ \gamma \ \textrm{is complex and} \ \delta\ \textrm{is real, with the condition that}\ \delta =0 \ \textrm{if} \ \gamma\in \mathbb R\ .\ \ \ \ \ \ \ \ \ \ \ \ \ \ \]

\vspace{0,2 cm}

\centerline{{\bf Main Result}}

\vspace{0,3 cm}

\hspace{-0,6 cm}{\bf Theorem 3.} The equation $P\big(\frac{d}{dt}\big)x = R(t)e^{\gamma t + i\delta}$ 
 has a particular solution $Q(t)e^{\gamma t + i\delta}$, where $Q(t)$ is an arbitrary polynomial satisfying, 
\[ (3.1)\ \ \ \ \ \ \ \ \ \ \ \frac{p^{(n)}(\gamma)}{n!}Q^{(n)} \ +\  ...\ +\   \frac{p'(\gamma)}{1!}Q' \ +\  \frac{p(\gamma)}{0!}Q = R.\ \ \ \ \ \ \ \ \ \ \ \ \ \ \ \ \ \ \ \ \ \ \ \ \ \ \ \ \ \ \ \ \ \  \]

Moreover,
\begin{itemize}

\item[(a)] If $\gamma \in \mathbb R$, then we can suppose $Q$ real. So, is real the solution $x(t)= Q(t)e^{\gamma t}$.
\item[(b)] If $\gamma \notin \mathbb R$, then $z(t)=Q(t)e^{\gamma t + i\delta}$ is a complex solution. If $\gamma = \alpha +\beta i$, where $\alpha,\beta \in \mathbb R$, then the functions $x(t)= \textrm{Re}[z(t)]$ and $y(t)= \textrm{Im}[z(t)]$ satisfy
\[ P\big(\frac{d}{dt}\big)x = R(t)e^{\alpha t}\cos(\beta t +\delta)\ \ \ ,\ \ \ P\big(\frac{d}{dt}\big)y = R(t)e^{\alpha t}\sin(\beta t +\delta)\  .\]
\item[(c)] If $p(\gamma) \neq 0$, then degree$(Q)=$ degree$(R)$.
\item[(d)] If $\gamma$ is a root of multiplicity $k$, then we can choose $Q(t) = t^kQ_1(t)$, with \,degree$(Q_1)=$ degree$(R)$.
\end{itemize}
\newpage
\hspace{-0,5 cm}{\bf Proof.} Let us apply Lemma 2. Searching for a solution  $ Q(t)e^{\gamma t +i\delta}$ we arrive at
\[ P\big(\frac{d}{dt}\big)\left\{Q(t)e^{\gamma t +i\delta}\right\} = \left[\frac{p^{(n)}(\gamma)}{n!}Q^{(n)}  +  ...+   \frac{p'(\gamma)}{1!}Q' +  p(\gamma)Q  \right]e^{\gamma t +i\delta}  =R(t)e^{\gamma t + i\delta} , \]
and so, we get equation (3.1). Hence, applying Lemma 1 and solving a trivial lower triangular linear system we find a polynomial solution $Q(t)$ of (3.1).

\vspace{0,2 cm}

Now, we proceed to verify the mentioned properties.

\begin{itemize}
\item[(a)] If $\gamma$ is real, it is obvious that we can suppose $Q$ also real.
\item[(b)] It is obvious.
\item[(c)] It is obvious. 
\item[(d)] If $\gamma$ is a root of multiplicity $k$ then the equation (3.1) becomes,
\[ \frac{p^{(n)}(\gamma)}{n!}Q^{(n)} \, +\,  ...\, +\, \frac{p^{(k)}(\gamma)}{k!}Q^{(k)}\ = R(t),\ \ \ p^{(k)}(\gamma) \neq 0\ ,\]
which has a polynomial solution $y=Q^{(k)}$, degree$(Q^{(k)})=$ degree($R$). Integrating $k$ times the function $y(t)$, but choosing at each time the independent term as zero, we find a suitable polynomial solution of (3.1) $\blacksquare$
\end{itemize}

\vspace{0,4 cm}

\hspace{-0,6 cm}{\bf Examples}

\begin{itemize}

\item[(E1)] Let us solve the ode $x'' - 2x' + 2x = t^2 e^t\sin (3t+5)$, where $x=x(t): \mathbb R \to \mathbb R$.

The characteristic polynomial $p(\lambda)$ satisfies,
\[p(\lambda)= \lambda^2 -2\lambda +2= (\lambda-1)^2 + 1\ ,\ p'(\lambda) = 2(\lambda -1)\ \textrm{and}\ p''(\lambda)=2\ .\]
 The general solution of the homogeneous associated equation is
\[  c_1e^t \cos t + c_2e^t\sin t\,,\ \ \ c_1\,,c_2 \in \mathbb R .\] 
By the theorem, the complex equation
\[ z'' -2z' + 2z\ =\  t^2e^{(1+3i)t +5i}\,,\ \  \    \]
has a solution $z(t)  =Q(t)e^{(1+3i)t +5i}$, 
with $Q$ a polynomial that satisfies
\[(E1.1)\ \ \ \ \ \ Q'' + p'(1+3i)Q' + p(1+3i)Q  = t^2\,.\ \ \ \ \ \ \ \ \ \ \ \ \ \ \ \ \ \ \ \ \ \ \ \ \ \ \ \ \ \ \ \ \ \ \ \ \ \]
Substituting $p(1+3i)=-8$ and $p'(1+3i)= 6i$ in equation (E1.1) we get,
\[ Q'' -8 Q' +6iQ=t^2\ .\]
So, we have that $Q(t) = \frac{t^2}{6i} + At +B$, with $A, B\in \mathbb C$. It is easy to see that,
\[Q(t) = -\frac{i}{6}t^2 -\frac{4}{9}t + \frac{1}{18} +\frac{16}{27}i\ .\]
Hence, we find that 
\begin{displaymath}
\left\{\begin{array}{ll}
z(t) &=  \Big[-\frac{i}{6}t^2 -\frac{4}{9}t + \frac{1}{18} +\frac{16}{27}i\Big]e^t\,[\cos(3t+5) + i\sin(3t+5)]\\
\\
\textrm{Im}[z(t)]&= e^t\Big[- \frac{t^2\cos(3t+5)}{6} \, -\,\frac{4t\sin (3t+5)}{9} \,  +\, \frac{\sin(3t+5)}{18}\, + \, \frac{16\cos(3t+5)}{27}\Big]\ 
\end{array}
\right.
\end{displaymath}

The general solution of the given ode is,
\[x(t) \ =\ c_1e^t \cos t + c_2e^t\sin t \ + \ \ \ \ \ \ \ \ \ \ \ \ \ \ \ \ \ \ \ \ \ \ \ \ \ \ \ \ \ \ \ \ \  \ \ \ \ \ \ \ \ \ \ \ \ \ \ \ \ \ \ \ \ \ \ \ \ \ \ \ \ \ \ \ \]
\[\ \ \ \ \ \ +\ e^t\Big[\Big(- \frac{t^2}{6} +\frac{16}{27}\Big)\cos(3t+5)  +\Big(-\frac{4t}{9}  + \frac{1}{18}\Big)\sin(3t+5) 
\Big]\,, \ \ \ \ \ \ \ \]
where $c_1$ and $c_2$ are arbitrary real constants.

\item[(E2)] Let us solve the ode $x''' - 5x'' + 3x' + 9x = t^5 e^{3t}$, where $x = x(t): \mathbb R \to \mathbb R$. 

The characteristic polynomial is 
\[p(\lambda) = \lambda^3 - 5\lambda^2 + 3\lambda + 9 = (\lambda - 3)^2(\lambda + 1)\ .\]
The general solution of the associated homogeneous equation is 
\[c_1e^{3t} + c_2te^{3t} + c_3e^{-t}\,,\ c_1, c_2,c_3 \in \mathbb R \ .\]
By the theorem above, the ode has a particular solution $Q(t)e^{3t}$ such that,
\[\frac{p'''(3)}{3!}Q''' + \frac{p''(3)}{2!}Q'' + \frac{p'(3)}{1!}Q' + \frac{p(3)}{0!}Q = t^5 \ .\]
Since $p' = 3\lambda^2 - 10\lambda + 3$, $p'' = 6\lambda - 10$ and $p''' = 6$, we easily get
\[Q''' + 4Q'' = t^5 \ .\]
Substituting $y = Q''$ in this last equation we find the ode $y' + 4y = t^5$, which has a solution $y= \frac{t^5}{4} + at^4 + bt^3 + ct^2 + dt + e$, for certains $a,b,c,d,e \in \mathbb R$. It is easy to verify that 
\[y(t)= \frac{t^5}{4} - \frac{5t^4}{16} + \frac{5t^3}{16} - \frac{15t^2}{64} + \frac{15t}{128} - \frac{15}{512}\ .\]
Hence, by integrating $y(t)$, we can choose
\[Q = \frac{t^7}{168} -\frac{t^6}{96} + \frac{t^5}{64} -\frac{5t^4}{256} + \frac{5t^3}{256} - \frac{15t^2}{1024}\ .\]
Therefore, the general solution of the given ode is,
\[x (t)= c_1e^{3t} + c_2te^{3t} + c_3e^{-t} + \Big(\frac{t^7}{168} -\frac{t^6}{96} + \frac{t^5}{64} -\frac{5t^4}{256} + \frac{5t^3}{256} - \frac{15t^2}{1024}\Big)e^{3t} \,,\  c_i's \in \mathbb R\, . \]

\end{itemize}

\

\hspace{- 0,6 cm}{\bf Acknowledgements.} I express my gratitude to Prof. Anthony W. Knapp for his suggestions and commentaries.

\

\centerline{{\bf References}}

\vspace{0,4 cm}

\hspace{- 0,6 cm}[1.] Apostol, T. M., {\sl Differential Equations and Their Applications}, Springer-Verlag, 1975.

\hspace{- 0,6 cm}[2.] Birkhoff, G. and Rota, G. C.,  {\sl Ordinary Differential Equations}, 4th ed.,  John Wiley \& Sons, 1989.

\hspace{- 0,6 cm}[3.] Braun, M., {\sl Differential Equations and Their Applications}, Springer-Verlag, 1975.

\hspace{- 0,6 cm}[4.] Coddington, E. A. and Levinson, N., {\sl Theory of Ordinary Differential Equations}, McGraw-Hill, 1955.

\hspace{- 0,6 cm}[5.] Krantz, S. G., and Simmons, G. F., {\sl Differential Equations}, The Walter Rudin Student Series in Advanced Mathematics, McGraw-Hill, 2007.

\hspace{-0,5 cm}[6.] Oliveira, O. R. B., ``The Fundamental Theorem of Algebra: An Elementary and Direct Proof'', {\sl{ The Mathematical Intelligencer}} 33, No. 2, (2011), 1-2.

\hspace{- 0,6 cm}[7.] Pontryagin, L., {\sl Ordinary Differential Equations}, Addison-Wesley, 1962.

\hspace{- 0,6 cm}[8.] Robinson, J. C., {\sl Ordinary Differential Equations}, Cambridge University Press, 2007.

\hspace{- 0,6 cm}[9.] Ross, Clay C., ``Why the Method of Undetermined Coefficients Works'', {\sl American Mathematical Monthly} 98 (1991), pp. 747-749.

\hspace{- 0,6 cm}[10.] Ross, Clay C., {\sl Differential Equations: an Introduction with Mathematica}, 2nd ed., Undergraduate Series In Mathematics, Springer, 2010.

\hspace{- 0,6 cm}[11.] Zill, D. G., {\sl A First Course in Differential Equations with Applications}, 4th ed., Prindle, Weber \& Schmidt, 1989.

\

\

\hspace{- 0,6 cm} Oswaldo Rio Branco de Oliveira

\hspace{- 0,5 cm}Departamento de Matemática 

\hspace{- 0,5 cm}Universidade de São Paulo 

\hspace{- 0,5 cm}Rua do Matão 1010, CEP 05508-090

\hspace{- 0,5 cm}São Paulo, SP

\hspace{- 0,5 cm}Brasil

\vspace{0,2 cm}

\hspace{- 0,5 cm}e-mail: oliveira@ime.usp.br

\end{document}